\documentclass{amsart}
\usepackage{amssymb,latexsym}
\usepackage[dvips]{graphicx}
\newtheorem{theorem}{Theorem}[section]
\newtheorem{prop}[theorem]{Proposition}
\newtheorem{lemma}[theorem]{Lemma}

\newtheorem{definition}[theorem]{Definition}
\newtheorem{cor}[theorem]{Corollary}

\def\co{\colon\thinspace}

\begin{document}

\title{Minimality and symplectic sums}\author{Michael Usher}
\address{Department of Mathematics\\ Princeton University\\ Princeton, NJ  08540}
\email{musher@math.princeton.edu} \subjclass{53D35, 53D45}

\begin{abstract} Let $X_1, X_2$ be symplectic $4$-manifolds
containing symplectic surfaces $F_1,F_2$ of identical positive genus
and opposite squares.  Let $Z$ denote the symplectic sum of $X_1$
and $X_2$ along the $F_i$.  Using relative Gromov--Witten theory, we
determine precisely when the symplectic $4$-manifold $Z$ is minimal
(\emph{i.e.}, cannot be blown down); in particular, we prove that
$Z$ is minimal unless either: one of the $X_i$ contains a
$(-1)$-sphere disjoint from $F_i$; or one of the $X_i$ admits a
ruling with $F_i$ as a section.  As special cases, this proves a
conjecture of Stipsicz asserting the minimality of fiber sums of
Lefschetz fibrations, and implies that the non-spin examples
constructed by Gompf in his study of the geography problem are
minimal.
\end{abstract}
 \maketitle
\section{Introduction}
Let $(X_1,\omega_1)$, $(X_2,\omega_2)$ be symplectic $4$-manifolds,
and let $F_1\subset X_1$, $F_2\subset X_2$ be two-dimensional
symplectic submanifolds with the same genus whose homology classes
satisfy $[F_1]^2+[F_2]^2=0$, with the $\omega_i$ normalized to give
equal area to the surfaces $F_i$. For $i=1,2$, a neighborhood of
$F_i$ is symplectically identified by Weinstein's symplectic
neighborhood theorem \cite{W} with the disc normal bundle $\nu_i$ of
$F_i$ in $X_i$.  Choose a smooth isomorphism $\phi$ of the normal
bundle to $F_1$ in $X_1$ (which is a complex line bundle) with the
dual of the normal bundle to $F_2$ in $X_2$. According to \cite{G}
(and independently \cite{MW}), the \emph{symplectic sum}
\[ Z=X_1\#_{F_1=F_2}X_2=(X_1\setminus
\nu_1)\cup_{\partial\nu_1\sim_{\phi}\partial\nu_2}(X_2\setminus
\nu_2)\] carries a natural deformation class of symplectic
structures (note that the diffeomorphism type of $Z$ may depend on
the identification $\phi$, as seen for instance in Example 3.2 of
\cite{G}, a feature which is mostly suppressed from the notation
hereinafter).

Symplectic sums along $S^2$ are well-understood; according to pp.
563--566 of \cite{G} such a symplectic sum amounts to either blowing
down a sphere of square $-1$ or $-4$ in one of the summands, taking
the fiber sum of two ruled surfaces, or leaving the diffeomorphism
type of one of the summands unchanged, and then possibly blowing up
the result. Accordingly, we shall restrict our attention to
symplectic sums along surfaces of positive genus.

Recall that a symplectic $4$-manifold $M$ is called \emph{minimal}
if it contains no symplectically embedded spheres of square $-1$,
and hence cannot be expressed as a blowup of another symplectic
$4$-manifold.  According to results arising from Seiberg--Witten
theory (\cite{T},\cite{Li2}), this is equivalent to the condition
that $M$ not contain any \emph{smoothly} embedded spheres of square
$-1$.

In this note, we resolve completely the question of under what
circumstances the above symplectic sum $Z$ is minimal.  Our result
may be summarized as follows:
\begin{theorem} \label{main} Let the symplectic sum $Z=X_1\#_{F_1=F_2} X_2$ be formed as
above, and assume that the $F_i$ have positive genus $g$.   Then:
\begin{itemize}
\item[(i)] If either $X_1\setminus F_1$ or $X_2\setminus F_2$
contains an embedded symplectic sphere of square $-1$, then $Z$ is
not minimal.
\item[(ii)] If one of the summands $X_i$ (for definiteness, say
$X_1$) admits the structure of an $S^2$-bundle over a surface of
genus $g$ such that $F_1$ is a section of this fiber bundle, then
$Z$ is minimal if and only if $X_2$ is minimal.
\item[(iii)] In all other cases, $Z$ is minimal.
\end{itemize}\end{theorem}

Case (i) above should be obvious: if $X_1$ admits a sphere of square
$-1$ which misses $F_1$, then cutting out a small neighborhood of
$F_1$ and replacing the neighborhood with something else will not
change that fact. In the situation of Case (ii), note that since
$F_1$ and $F_2$ have opposite squares, the complement of a
neighborhood of $F_1$ in the ruled surface $X_1$ is diffeomorphic to
a neighborhood of $F_2$, and so effectively the symplectic sum
operation cuts out a neighborhood of $F_2$ and then glues it back in
via a map which (since $\phi\co \partial \nu_1\to \partial \nu_2$ is
a bundle isomorphism, not just a diffeomorphism) takes a meridian of
$F_2$ to itself; certainly for some choices of the gluing map the
result is just diffeomorphic to $X_2$ and so of course will have the
same minimality properties, and indeed the author is unaware of any
cases in which the diffeomorphism type of $X_2$ is changed by
summing with a ruled surface along a section with a nonstandard
 choice of gluing map.

Case (iii) is thus quite broad, and in particular confirms the
belief expressed by the authors of \cite{LS} that the main theorem
of that paper could be generalized. For some prior results asserting
the minimality of symplectic sums in various special cases, see
\cite{S1}, Theorem 1.5 of \cite{S2}, Theorem 2.5 of \cite{P},
\cite{LS}, and \cite{Li}.

We now turn to some corollaries of Theorem \ref{main}.  Recall that
a map $f\co X\to \Sigma$ from a symplectic $4$-manifold to a
$2$-manifold is called a symplectic Lefschetz fibration provided
that its smooth fibers are symplectic submanifolds and $f$ has just
finitely many critical points, near each of which it is given in
orientation preserving complex coordinates by $(z_1,z_2)\mapsto
z_1z_2$.  $f$ is called relatively minimal if none of its singular
fibers (each of which is a nodal curve) contains a $(-1)$-sphere as
a reducible component.  Note that in general a Lefschetz fibration
admits compatible symplectic structures provided that its fiber is
homologically essential; this condition is automatic if the fiber
genus is at least two or if the genus is one and the map has at
least one critical point.  The fiber sum of two symplectic Lefschetz
fibrations whose fibers have the same genus is just the symplectic
sum along a smooth fiber.  The following was conjectured by A.
Stipsicz in \cite{S1}:

\begin{cor} \label{lef} Let $f_i\co X_i\to S^2$ ($i=1,2$) be
relatively minimal symplectic Lefschetz fibrations on $4$-manifolds
$X_1,X_2$ whose fibers $F_1,F_2$ have the same positive genus $g$,
and assume that neither $f_i$ is the projection $\Sigma_g\times
S^2\to S^2$. Then the fiber sum $X_1\#_{F_1=F_2} X_2$ is
minimal.\end{cor}

 Another consequence of Theorem \ref{main} relates to the geography
of minimal symplectic $4$-manifolds, that is, the question of which
pairs of integers $(a,b)$ have the property that there is a minimal
symplectic $4$-manifold $M_{a,b}$ such that $c_{1}^{2}(M_{a,b})=a$
and $c_2(M_{a,b})=b$.  Any such pair $(a,b)$ necessarily satisfies
$a+b\equiv 0\, (mod\, 12)$ by the Noether formula.  By performing
symplectic sums on certain manifolds covered by Case (iii) of
Theorem \ref{main}, R. Gompf in \cite{G} realized a great many such
pairs as the Chern numbers of symplectic manifolds with prescribed
fundamental group; in the case that the Chern numbers are consistent
with Rohlin's theorem (so that $(a,b)$ has form $(8k, 4k+24l)$)
Gompf was able to arrange that the resulting manifold be spin and
hence minimal, but in other cases minimality appeared likely but
could not be proven. However, Theorem \ref{main} in conjunction with
Remark 2 after Theorem 6.2 of \cite{G} allow us to deduce:
\begin{cor} \label{Gompfex}Let $G$ be any finitely presentable group. There is a
constant $r(G)$ with the property that if $(a,b)\in \mathbb{Z}^2$
satisfies $a+b\equiv 0\, (mod \,12)$ and $0\leq a\leq 2(b-r(G))$
then there is a minimal symplectic $4$-manifold $M_{a,b,G}$ such
that $\pi_1(M_{a,b,G})\cong G$, $c_{1}^{2}(M_{a,b,G})=a$, and
$c_2(M_{a,b,G})=b$.\end{cor}

 Note that most of Gompf's examples involve taking symplectic sums
in which at least one of the summands is rational, so these examples
are not covered by previous results on the minimality of symplectic
sums.  In the simply connected case, different constructions have
been used to obtain symplectic manifolds occupying large parts of
the $c^{2}_{1}c_2$-plane and then to show that they are minimal
using gauge theory; see, \emph{e.g.}, Theorem 10.2.14 of \cite{GS}.

The next two sections are occupied with the proof of Theorem
\ref{main}. This proof splits into two parts: first, we use relative
Gromov--Witten theory to give a condition on the pairs
$(X_1,F_1),(X_2,F_2)$ in terms of the intersection numbers of the
$F_i$ with  holomorphic spheres which is sufficient to guarantee the
minimality of $Z$ (namely, the $F_i$ should be ``rationally
$K$-nef,'' defined below). We then see that surfaces of positive
genus in symplectic $4$-manifolds are always rationally $K$-nef
except in the cases (i) and (ii) in Theorem \ref{main}; this follows
from results of Seiberg--Witten theory concerning the canonical
class (\cite{T},\cite{Liu}) when the ambient manifold is not a
blowup of a rational or ruled surface, while the ruled case can be
handled fairly directly and the rational case depends in part on the
analysis of the chamber structure in the cohomology of rational
surfaces that was carried out in \cite{FM}.  Finally, in the last
section we prove Corollaries \ref{lef} and \ref{Gompfex}.

\subsection*{Acknowledgements} I am indebted to T.J. Li for finding
an oversight in an earlier version of this paper.  This work was
partially supported by an NSF Postdoctoral Fellowship.

\section{Symplectic sums along rationally K-nef surfaces are
minimal}\label{relgr}

If $(X,\omega)$ is a symplectic $4$-manifold, $g\geq 0$, $A\in
H_2(X;\mathbb{Z})$, and $\alpha\in H^*(X;\mathbb{Z})$, we let
$GW_{g,A}^{X}(\alpha)$ denote the Gromov--Witten invariant
(\cite{RT}) counting perturbed-pseudoholomorphic maps from a surface
of genus $g$ into $X$, representing the homology class $A$, and
passing through a cycle Poincar\'e dual to $\alpha$.

Suppose now that the symplectic sum $Z=X_1\#_{F_1=F_2}X_2$ is
\emph{not} minimal, so that $Z$ contains a symplectic sphere $E$ of
square $-1$. Using an almost complex structure $J$ generic among
those making $E$ pseudoholomorphic, we immediately see that
$GW_{0,[E]}^{Z}(1)=\pm 1$, since the operator
$\overline{\partial}_{J}$ will (for generic $J$) have nondegenerate
linearization at the embedding of $E$, and positivity of
intersections (see, \emph{e.g.}, Theorem E.1.5 of \cite{MS})
prevents the existence of any other $J$-holomorphic spheres
homologous to $E$. We now review how to use this nonvanishing to
deduce the existence of nonzero relative Gromov--Witten invariants
in certain homology classes in $X_1, X_2$ by means of the gluing
results in \cite{IP}.

The main theorem of \cite{IP} (or, similarly, that of \cite{LR})
provides a somewhat complicated formula expressing the
Gromov--Witten invariants of $Z$ in terms of the Gromov--Witten
invariants of $X_1$ and $X_2$ \emph{relative} to the fibers $F_1$
and $F_2$, the latter invariants having been defined in \cite{IP1}.
Without reproducing the formula, we recall the essential points,
referring readers to \cite{IP} for details: in forming the
symplectic sum, letting $x_1$ and $x_2$ be complex coordinates in
the normal bundles $\nu_1,\nu_2$, one performs the identification by
taking $x_1x_2=\lambda$ for some small $\lambda\in \mathbb{C}^*$;
this results in a symplectic $6$-manifold $\mathcal{Z}$ equipped
with a projection to the disc whose fiber
$(Z_{\lambda},\omega_{\lambda})$ over $\lambda\in D^2\setminus\{0\}$
is isotopic to $Z$ but whose fiber over $0$ is
$Z_0=X_1\cup_{F_1=F_2} X_2$. As the parameter $\lambda$ approaches
zero, pseudoholomorphic curves in $Z_{\lambda}$ limit to trees of
curves in $Z_0$ consisting of: curves in $X_1$ meeting $F_1$ in
isolated points; curves contained in the identified surfaces
$F=F_1=F_2$; and curves in $X_2$ meeting $F_2$ in isolated points.
In fact there is a quantifiable finite-to-one correspondence between
these trees (which are counted by a combination of (relative)
Gromov--Witten invariants in $X_1, F,$ and $X_2$) and the curves
counted by the Gromov--Witten invariants of $Z$.  This leads to a
gluing formula (Theorem 12.4 of \cite{IP}) expressing the latter
invariants in terms of the former.

Let us return to our case, where $Z$ admits a sphere $E$ of square
$-1$, and consider what the gluing result of \cite{IP} allows us to
deduce about the Gromov--Witten invariants of $X_1$ and $X_2$.
First, notice that since $E$ is a sphere the limiting tree discussed
in the previous paragraph will consist only of genus-zero curves,
and so will not have any components mapped into $F$, since
$\pi_2(F)=0$.  So the gluing formula will express the nonzero
Gromov--Witten invariant $GW_{0,[E]}^{Z}(1)$ in terms of invariants
which count pairs $(C_1,C_2)$ of possibly-disconnected curves (each
of whose components are spheres) representing elements of the
intersection-homology spaces $\mathcal{H}_{X_1}^{F_1}$,
$\mathcal{H}_{X_2}^{F_2}$ defined in Section 5 of \cite{IP1}, and
which have matching intersections with the fiber $F=F_1=F_2$. So we
obtain nonzero relative invariants in both $X_1$ and $X_2$, which in
particular count (generally disconnected) holomorphic curves in
certain classes $A_i\in H_2(X_i;\mathbb{Z})$ $(i=1,2)$ (subject to
some additional constraints on their intersections with $F_i$).  We
note also that while in the general situation considered in
\cite{IP} the Gromov--Witten invariants are counts of maps
satisfying a perturbed Cauchy-Riemann equation
$\overline{\partial}_Ju=\nu$, in our case we can take $\nu=0$ by
virtue of the fact that $[E]$ is a primitive class in
$H_2(Z;\mathbb{Z})$ (so that the relevant moduli space has no strata
consisting of multiple covers). Thus the curves $C_1$ and $C_2$ will
be genuinely pseudoholomorphic curves for almost complex structures
$J_1,J_2$ on $X_1,X_2$, which may be chosen generically among those
pairs of almost complex structures on the $X_i$ which preserve
$TF_i$ and agree on the identified neighborhoods $\nu_i$ of $F_i$.
Incidentally, in our proof of Theorem \ref{main} we will in fact
only need the fact that a nonvanishing invariant $GW_{0,[E]}^{Z}(1)$
gives rise rise via a Gromov-type compactness theorem to (generally
reducible, non-reduced) $J_i$-holomorphic curves $C_i$ for some
almost complex structures $J_i$ which make the $F_i$
pseudoholomorphic.  The full strength of the Ionel-Parker theorem,
which shows that the $C_i$ are in fact enumerated by nonvanishing
relative invariants, is not needed for this conclusion.

Let $d=A_1\cap[F_1]=A_2\cap[F_2]$.  Then according to Lemma 2.2 of
\cite{IP}, we have, where for a symplectic manifold $M$ we denote
the canonical class of $M$ by $\kappa_M$,
\[ \langle\kappa_Z,[E]\rangle=\langle\kappa_{X_1},
A_1\rangle+\langle\kappa_{X_2},A_2\rangle+2d.\]  But $E$ is an
embedded $(-1)$ sphere and so satisfies
$\langle\kappa_Z,[E]\rangle=-1$ by the adjunction formula; thus we
may suggestively rewrite the above equation as
\begin{equation} \label{can}
\langle\kappa_{X_1}+PD[F_1],A_1\rangle+\langle\kappa_{X_2}+PD[F_2],A_2\rangle=-1.\end{equation}

Accordingly, we make the following definition (recall from \cite{MS}
that a pseudoholomorphic curve $u\co\Sigma \to X$ is called simple
if no two disjoint open sets in its domain $\Sigma$ have the same
image; in this case the map $u$ is generically injective):

\begin{definition} \label{knef} Let $(X,\omega)$ be a symplectic four-manifold.  An embedded symplectic surface $F\subset X$ is called
rationally $K$-nef if, whenever $J$ is an almost complex structure
preserving $TF$ and $A\in H_2(X;\mathbb{Z})$ is represented by a
simple $J$-holomorphic sphere, we have
\begin{equation} \label{defnef} \langle \kappa_X+PD[F],A\rangle\geq
0,\end{equation} where $\kappa_X$ is the canonical class of
$X$.\end{definition} The classes  $A_1,A_2$ in (\ref{can}) are both
represented by $J_i$-holomorphic curves for appropriate almost
complex structures $J_i$ which make the respective $F_i$ holomorphic
(indeed, the $A_i$ are the homology classes underlying elements of
$\mathcal{H}^{F_i}_{X_i}$ having nonvanishing (possibly
disconnected) relative Gromov--Witten invariants), and  each
component of these pseudoholomorphic curves has genus zero and is
either simple or a multiple cover of a simple curve. So if $F_1$ and
$F_2$ are both rationally $K$-nef then (\ref{can}) cannot hold.
Thus:

\begin{prop} \label{knefmin} If $F_1\subset X_1, F_2\subset X_2$
are rationally $K$-nef surfaces, then the symplectic sum
$X_1\#_{F_1=F_2} X_2$ is minimal.\end{prop}

\section{Many surfaces are rationally K-nef}

To prove Case (iii) of our main theorem, we need to demonstrate
that, if $F\subset X$ is an embedded symplectic surface of positive
genus which intersects all symplectically embedded $(-1)$-spheres in
$X$, and which
 is not a section of
a ruled surface, then $F$ is rationally $K$-nef. Let us begin with
some general observations.  We wish to show that if $J$ preserves
$TF$ and $A\in H_2(Z;\mathbb{Z})$ is represented by a simple
$J$-holomorphic sphere, then $\langle \kappa_X+PD(F),A\rangle \geq
0$.  Now by positivity of intersections $[F]\cap A\geq 0$, so in the
case that $\langle \kappa_X, A\rangle\geq 0$ the desired inequality
obviously holds. Further, since $F$ meets all embedded symplectic
spheres of square $-1$, if $A$ is represented by such a sphere then
positivity of intersections implies that $[F]\cap A\geq 1$, and so
since in this case $\langle \kappa_X,A\rangle =-1$ the inequality
holds when $A$ is represented by an embedded $(-1)$-sphere as well.

Now in general the adjunction formula (\emph{e.g.}, Corollary E.1.7
of \cite{MS}) shows that $A^2+\langle\kappa_X,A\rangle\geq -2$, with
equality if and only if $A$ is represented by an embedded sphere.
From this we see that $A^2\geq 0$ unless either $\langle
\kappa_X,A\rangle\geq 0$ or $A^2=\langle\kappa_X,A\rangle=-1$ (in
which case $A$ is represented by an embedded $(-1)$-sphere), and the
latter two cases have already been dispensed with.  Hence:

\begin{lemma}\label{adj} An embedded symplectic surface $F\subset X$ which intersects all symplectically embedded $(-1)$-spheres in $X$ is
rationally $K$-nef provided that, whenever $J$ is an almost complex
structure preserving $TF$ and $A\in H_2(X;\mathbb{Z})$ is
represented by a simple $J$-holomorphic sphere such that \[ A^2\geq
0 \quad \mbox{and} \quad \langle\kappa_X,A\rangle <0,\]  we have
\[ \langle \kappa_X+PD[F],A\rangle\geq 0.\] \end{lemma}

The task of showing that the surfaces $F$ in question are rationally
$K$-nef now naturally splits up into several cases.

\subsection{Case 1: $b^+(X)> 1$} In this case, according to results
of \cite{T}, the Gromov--Taubes invariant of the canonical class is
nonzero, and so for generic almost complex structures among those
preserving $TF$ the class $PD(\kappa_X)$ is represented by a
(possibly disconnected) pseudoholomorphic curve, whose only
spherical components are embedded $(-1)$-spheres.  For an arbitrary
complex structure $J$ preserving $TF$, then, Gromov compactness
shows that $PD(\kappa_X)$ is represented at least by a union of
pseudoholomorphic bubble trees.  So if $A$ is represented by a
simple $J$-holomorphic sphere with nonnegative square, then
positivity of intersections between the representatives of
$PD(\kappa_X)$ and $A$ shows that $\langle \kappa_X,A\rangle\geq 0$
(of course, if $A$ had negative square the representative of
$PD(\kappa_X)$ might contain the representative of $A$ as a
component, allowing the intersection to be negative; indeed when $A$
is the class of an embedded $(-1)$-sphere this is precisely what
happens). So the sufficient condition provided by Lemma \ref{adj} is
vacuously satisfied and:
\begin{prop} If $b^+(X)>1$ and $F$ is an embedded symplectic surface
in $X$ which intersects all symplectically embedded $(-1)$-spheres
in $X$, then $F$ is rationally $K$-nef.\end{prop}

\subsection{Case 2: $X$ is a (possibly trivial) blowup of an irrational ruled
surface}   First suppose $X=(\Sigma\times
S^2)\#n\overline{\mathbb{C}P^2}$, where $\Sigma$ has genus $h>0$.
Let $\sigma$ denote the homology class of the strict transform of
$\Sigma\times\{pt\}\subset \Sigma\times S^2$, $f$ the homology class
of the strict transform of $\{pt\}\times S^2$, and $e_1,\ldots,e_n$
the classes of the $n$ exceptional divisors.  We then have \[
PD(\kappa_X)=-2\sigma+(2h-2)f+\sum_{i=1}^{n}e_i. \]  Now assuming
that $n>0$ the $(-1)$-spheres $e_i$ all have nonvanishing
Gromov--Witten invariant, so if $F$ meets all embedded
$(-1)$-spheres we have $[F]\cap e_i\geq 1$ for each $i$; further  we
get $(-1)$-spheres in each of the classes $f-e_i$ by choosing a
symplectic sphere homologous to $\{pt\}\times S^2$ in $\Sigma\times
S^2$ which meets the $i$th blow-up point and none of the others, and
taking its strict transform, so that $[F]\cap (f-e_i)\geq 1$  as
well. Also, regardless of whether $n$ is positive, the class $f$ has
a nonvanishing Gromov--Witten invariant with a single point
constraint, and so by considering an almost complex structure
preserving $TF$ we see that by positivity of intersections $[F]\cap
f\geq 0$ (since $F$ has positive genus, no rational curve can share
a component with it).  Hence \[
[F]=c\sigma+df-\sum_{i=1}^{n}a_ie_i\] where $c\geq 0$ and, for each
$i$, $1\leq a_i<c$. Since $h>0$, the only rational curves which
appear  will be contained in the fiber of the ruling and so, if
their Chern number is positive, will have homology class in the cone
spanned by the classes $e_i$ and $f-e_i$, if $n>0$, or will be a
multiple of $f$ if $n=0$. If $n>0$, then since $e_i$ and $f-e_i$
have square $-1$ and so pair positively with $F$ and as $-1$ with
$\kappa_X$, $\kappa_X+PD[F]$ is nonnegative on each of them and
hence also on any element of the cone spanned by them.  So if $n>0$
$F$ is rationally $K$-nef.  In case $n=0$,  to verify that $F$ is
rationally $K$-nef we just need to ensure that $\kappa_X+PD[F]$
pairs positively with $f$, which amounts to the statement that
$c\geq 2$. Now we have $[F]=c\sigma+df\in H_2(\Sigma\times
S^2;\mathbb{Z})$ with $c\geq 0$; observe that if $c=0$, $F$ would be
homologous to a multiple of the fiber of the projection
$\Sigma\times S^2\to \Sigma$, and this multiple would be positive
since both $F$ and the fiber are symplectic.  So if $c=0$ we would
have $[F]=df$ with $d>0$, and so $2g(F)-2=-2d<0$, contradicting the
fact that $F$ has positive genus. If $c=1$, we see that \[
2g(F)-2=[F]^2+\langle\kappa_X,[F]\rangle=2d-2d+2h-2,\] \emph{i.e.},
$h=g(F)$ and $X$ is the total space of a (trivial) fibration over a
surface of genus $g$.  Thus $F$ is rationally $K$-nef unless $X$ is
a trivial $S^2$-fibration over a surface of genus $g=g(F)$ whose
fiber has homological intersection number $1$ with $F$.

Now suppose that $X$ is a nontrivial sphere bundle
$\Sigma\tilde{\times}S^2$ over a surface of genus $h>0$.  In this
case, the homology is generated by sections $s^+,s^-\in
H_2(X;\mathbb{Z})$ such that $(s^{\pm})^2=\pm 1$ and $s^+\cap
s^-=0$.  The fiber of the fibration represents the homology class
$f=s^+-s^-$, and the canonical class is given by
$PD(\kappa_X)=(2h-2)f-(s^++s^-)=(2h-3)s^+-(2h-1)s^-$. Again, since
$h>0$, the only classes represented by $J$-holomorphic rational
curves are multiples of the fiber class $s^+-s^-$; as such,
positivity of intersections implies that our surface $F$ represents
a class of form $[F]=cs^++ds^-$ where $c+d\geq 0$.  The desired
rational $K$-nef condition for $[F]$ is that $c+d\geq 2$, so we just
need to rule out $c+d\in\{0,1\}$.  We find
\[ 2g(F)-2=[F]^2+\langle
\kappa_X,[F]\rangle=(c+d)(c-d)+(2h-2)(c+d)+(d-c).\]  If $c+d=0$,
then for some $a\in\mathbb{Z}$, $[F]=a(s^+-s^-)=af$, with again
$a>0$ since $F$ is symplectic; then the adjunction formula would
give $2g(F)-2=-2a<0$, a contradiction.  If $c+d=1$, the adjunction
formula reads $2g(F)-2=2h-2$, so as before $g(F)=h$ and $F$ has
intersection number $1$ with the fibers.

We now make the following observation:
\begin{prop}\label{ruled}
 Suppose that $(X,\omega)$ is the total space of a ruled surface
over a surface of genus $g>0$, and that $F\subset X$ is a symplectic
submanifold having genus $g$ and homological intersection number one
with the fibers of the ruling.  Then $X$ admits a (possibly
different) ruling $\pi\co X\to \Sigma$ over a surface of genus $g$
such that $F$ is a section of $\pi$.\end{prop}
\begin{proof}This essentially follows from the techniques of
\cite{M}: let $J$ be an $\omega$-compatible almost complex structure
on $X$ which preserves $TF$.  The mere existence of an embedded
symplectic sphere of square zero in $X$ ensures that, letting
$\Sigma$ denote the moduli space of (unparametrized) $J$-holomorphic
spheres homologous to the fiber of the (original) ruling, $\Sigma$
is two-real dimensional, and the map $\pi\co X\to \Sigma$ which
takes a point $x\in X$ to the point of $\Sigma$ representing the
unique $J$-holomorphic sphere on which $x$ lies is an $S^2$-bundle
with symplectic (indeed $J$-holomorphic) fibers.  Since $F$ is also
$J$-holomorphic, the assumption on the intersection number ensures
that each of these fibers meets $F$ transversely and just once, so
letting $s\co \Sigma\to F$ denote the map which takes a $J$-curve
$C\in \Sigma$ to the unique point of $C\cap F$, we see that $s$ is a
section of $\pi$ with image $F$.  Of course, consideration of
$b_1(X)$ reveals that $\Sigma$, like $F$, has genus $g$.\end{proof}

Hence since, for $n>0$,
$(\Sigma\tilde{\times}S^2)\#n\overline{\mathbb{C}P^2}$ is
symplectomorphic to $(\Sigma\times S^2)\#n\overline{\mathbb{C}P^2}$,
we deduce:
\begin{prop} If the minimal model of $X$ is an irrational ruled
surface, and if $F\subset X$ is an embedded surface of positive
genus which intersects each embedded $(-1)$-sphere, then $F$ is
rationally $K$-nef unless $X$ is already minimal and admits a ruling
with $F$ as a section.\end{prop}

\subsection{Case 3: $X=S^2\times S^2$} Let  $[F]=c\sigma+df$ where $\sigma=[S^2\times\{pt\}]$ and
$f=[\{pt\}\times S^2]$. Then $PD(\kappa_X)=-2\sigma-2f$.
 Now there are nonzero Gromov--Witten
invariants counting holomorphic spheres in both the classes $\sigma$
and $f$, so since $F$ is symplectic by considering an almost complex
structure preserving $TF$ we see that $c,d\geq 0$ by positivity of
intersections between $F$ and the holomorphic spheres representing
$f$ and $\sigma$; of course $F$ is homologically essential, so $c$
and $d$ cannot both be zero. By the adjunction formula,
\[ 0\leq
2g(F)-2=[F]^2+\langle\kappa_X,[F]\rangle=2cd-2(c+d)=2((c-1)(d-1)-1),\]
and this then forces $c,d\geq 2$.  This implies the desired property
for $[F]=c\sigma+df$, since if $a\sigma+bf$ were some other class
represented by a $J$-holomorphic curve for some $J$ preserving $TF$,
positivity of intersections with the $J$-holomorphic representatives
of $f$ and $\sigma$ would imply that $a,b\geq 0$, and so $\langle
\kappa_X+PD[F],a\sigma+bf\rangle=a(c-2)+b(d-2)\geq 0$.  We have
shown:
\begin{prop} If $X=S^2\times S^2$ and $F\subset X$ is an embedded
symplectic surface of positive genus then $F$ is rationally
$K$-nef.\end{prop}

Before proceeding to the remaining cases, we recall a basic fact
about the intersection forms of symplectic $4$-manifolds
$(X,\omega)$ with $b^+(X)=1$. Where $n=b^-(X)$, the intersection
form on $H^2(X;\mathbb{Q})$ has type $(1,n)$. As such, the
``positive cone'' $\{\beta\in H^2(X;\mathbb{Q})|\beta^2>0\}$ has two
connected components.  A consequence of the Cauchy-Schwarz
inequality often called the ``light cone lemma'' then asserts that
the product of any two elements lying in the closure of the same
component of the positive cone is nonnegative, and in fact is
positive unless the elements are proportional and both have square
zero. In particular, one of the components (called the ``forward
positive cone'') is characterized by the property that all of its
elements pair positively with $[\omega]$, and so we obtain the
following fact which we shall make use of on two occasions:
\begin{lemma}\label{lightcone} If $b^+(X) = 1$ and $\alpha,\beta\in H^2(X;\mathbb{Q})$
satisfy $\alpha^2\geq 0$, $\beta^2\geq 0$, $\alpha\cdot [\omega]\geq
0$, and $\beta\cdot[\omega]\geq 0$, then $\alpha\cdot\beta\geq
0$.\end{lemma}

\subsection{Case 4: $X=\mathbb{C}P^2\#n\overline{\mathbb{C}P^2}$}

Let $H\in H_2(X;\mathbb{Z})$ be the class of the preimage of a
generic line in $\mathbb{C}P^2$ under the blowdown map
$X\to\mathbb{C}P^2$, and let $E_1,\ldots,E_n\in H_2(X;\mathbb{Z})$
be the homology classes of the exceptional divisors of the blowups.

\begin{lemma}\label{possquare} With the notation as above, let $\lambda\in
H^2(X;\mathbb{Z})$ satisfy $\lambda^2\geq \kappa_X\cdot\lambda$,
$\langle \lambda,H\rangle\geq 0$ and $\langle \lambda,[S]\rangle\geq
0$ for every embedded symplectic $(-1)$-sphere $S$.  Then
\[\lambda^2\geq 0.\]\end{lemma}
\begin{proof} If $n=0$, the intersection form is positive definite
and so the result is trivial.  If $n=1$, write
$\lambda=aPD(H)-bPD(E_1)$.  That
$\langle\lambda,H\rangle,\langle\lambda,E_1\rangle\geq 0$ shows that
$a,b\geq 0$, and that $\lambda^2\geq \kappa_X\cdot\lambda$ shows
that $a^2+3a\geq b^2+b$.  We wish to show that $a\geq b$; write
$b=a+k$, so we need to show that $k\leq 0$. Expanding out and
rearranging $a^2+3a\geq b^2+b$ gives \[ 2a\geq 2ak+k+k^2;\] noting
that $k\in\mathbb{Z}$ if $k$ were positive we would have $k\geq 1$
and so $2ak+k+k^2\geq 2a+2$, a contradiction. So indeed $k\leq 0$
and so $\lambda^2\geq 0$.

Assume for the rest of the proof that $n\geq 2$.  We hereinafter
view $\lambda$ as a real cohomology class rather than an integral
one. Of course, if $\kappa_X\cdot\lambda\geq 0$ then the hypothesis
of the lemma immediately shows $\lambda^2\geq 0$, so we may assume
that
\[ -\kappa_X\cdot \lambda >0.\]

Give $X$ an integrable complex structure $J$ which makes it a good
generic surface in the sense of \cite{FM}, \emph{i.e.} such that
$-\kappa_X$ is Poincar\'e dual to a smooth elliptic curve and all
 rational curves of negative square on $X$ are smooth and have square $-1$. (Such a $J$
exists by results from section I.2 of \cite{FM}.)   Let
$\mathcal{E}$ denote the set of all classes in $H_2(X;\mathbb{Z})$
which are represented by smooth rational curves  of square $-1$; of
course these spheres are symplectic, so we have
\[ \lambda\in \mathcal{S}:=\{\alpha\in
H^2(X;\mathbb{R})|-\kappa_X\cdot \alpha\geq
0,\langle\alpha,H\rangle\geq 0,\mbox{ and } (\forall E\in
\mathcal{E})(\langle\alpha,E\rangle\geq 0)\},\] and so the lemma
will be proven if we can show that the above set $\mathcal{S}$
consists only of elements of nonnegative square, which we now set
about doing.

Following section II.3 of \cite{FM}, let $\mathbf{H}(X)\subset
H_2(X;\mathbb{R})$ denote the set of elements of square $+1$, and
let $\mathcal{K}(X,J)\subset \mathbf{H}(X)$ denote the closure (in
$\mathbf{H}(X)$) of the set of cohomology classes of square 1 which
are represented by K\"ahler forms.  In \cite{FM}, $\mathcal{K}(X,J)$
is given two descriptions which are relevant to us.  On the one
hand, according to Proposition II.3.4 of that paper translated into
our notation, we have \[ \mathcal{K}(X,J)=\{\alpha\in
\mathcal{S}|\alpha^2=1\}.\]  On the other hand, Proposition II.3.6
shows that there is a \emph{chamber} $C_0(X)$, \emph{i.e.} the
closure of a connected component of \[\mathbf{H}(X)\setminus
\bigcup_{\{\alpha\in
H^2(X;\mathbb{Z}):\alpha^2=-1\}}\alpha^{\perp},\] with the property
that
\[ \mathcal{K}(X,J)\subset C_0(X)\] (in fact, $\mathcal{K}(X,J)$ has
a more explicit description as a ``$P$-cell'' of $C_0(X)$, but this
will not be important to us). We hence have $\{\alpha\in
\mathcal{S}|\alpha^2>0\}\subset \mathbb{R}_+\cdot C_0(X)$, whence
(since $\mathcal{S}$ is convex, and so if $\alpha\in \mathcal{S}$
has $\alpha^2=0$ a line segment from $\alpha$ to an element of
$\mathcal{S}$ with positive square gives a sequence of elements of
$\mathbb{R}_+\cdot C_0(X)$ approximating $\alpha$) \[
\{\alpha\in\mathcal{S}|\alpha^2\geq 0\}\subset
\overline{\mathbb{R}_+\cdot C_0(X)}.\]  We may now appeal to
Corollary II.1.20 of \cite{FM}, which asserts that, \emph{because
$n\geq 2$}, the set $\{\alpha\in \overline{\mathbb{R}_+\cdot
C_0(X)}|\alpha^2=0\}$ (and hence also $\{\alpha\in
\mathcal{S}|\alpha^2=0\}$) has no interior in $\{\alpha\in
H^2(X;\mathbb{R})|\alpha^2=0\}$.  So arguing just as in the proof of
Corollary II.1.21 of \cite{FM}, if $\mathcal{S}$ contained some
element $\beta$ of negative square, then since the intersection of
$\mathcal{S}$ with $\{\alpha^2>0\}$ (namely
$\mathbb{R}_+\cdot\mathcal{K}(X,J)$) has nonempty interior, the set
of points of square zero lying on line segments from $\beta$ to
elements of $\mathbb{R}_+\cdot\mathcal{K}(X,J)$ would sweep out a
set with nonempty interior in $\{\alpha^2=0\}$, which gives a
contradiction since $\mathcal{S}$ is convex.  Thus indeed all
elements of $\mathcal{S}$ have nonnegative square.\end{proof}

We now assume that $F$ is an embedded symplectic surface of positive
genus which intersects all embedded symplectic $(-1)$-spheres.
Write $[F]=aH-\sum_{i=1}^{n}b_iE_i$, so positivity of intersections
with the holomorphic representatives of $E_i$ for an almost complex
structure preserving $TF$ shows $b_i\geq 1$, and the fact that $F$
has positive area shows that $a>0$.  By the adjunction formula, we
have \[ 0\leq [F]^2+\langle\kappa_X,[F]\rangle=a(a-3)+\sum
(b_i-b_{i}^{2})\leq a(a-3),\] and so since $a>0$ we in fact have
$a\geq 3$.

Set \[ \lambda=\kappa_X+PD[F]=(a-3)PD(H)+\sum_{i=1}^{n}(1-b_i)E_i\in
H^2(X;\mathbb{Z})\]  So $\langle \lambda, H\rangle=a-3\geq 0$, and
if $S$ is an embedded symplectic $(-1)$-sphere, then $\langle
\lambda,[S]\rangle=[F]\cap[S]-1\geq 0$.  Further by the adjunction
formula $\lambda^2-\kappa_X\cdot\lambda=[F]^2+\langle
\kappa_X,[F]\rangle\geq 0$.  Hence by Lemma \ref{possquare}
$\lambda^2\geq 0$.

Assume that $A\in H_2(X;\mathbb{Z})$ is represented by a simple
$J$-holomorphic sphere for some almost complex structure $J$
preserving $TF$; in accordance with Lemma \ref{adj}, we assume that
$A^2\geq 0$.  The class $H$ has a nonvanishing Gromov--Witten
invariant with two point constraints, so positivity of intersections
guarantees that $A\cap H\geq 0$. So since $(\kappa_X+PD[F])^2\geq
0$, $A^2\geq 0$, and $\langle \kappa_X+PD[F],H\rangle$ and $A\cap H$
are both nonnegative, it follows from Lemma \ref{lightcone} that
also $\langle \kappa_X+PD[F],A\rangle \geq 0$.  Thus:
\begin{prop} If $F\subset \mathbb{C}P^2\#n\overline{\mathbb{C}P^2}$
is an embedded symplectic surface of positive genus $g$ which
intersects all embedded $(-1)$-spheres, then $F$ is rationally
$K$-nef.\end{prop}

Of course, since $(S^2\times
S^2)\#n\overline{\mathbb{C}P^2}=\mathbb{C}P^2\#(n+1)\overline{\mathbb{C}P^2}$
for $n\geq 1$, we've now exhausted all rational cases.

\subsection{Case 5: $b^+(X)=1$ and the minimal model of $X$ is neither rational nor ruled} We claim that no  $A\in
H_2(X;\mathbb{Z})$ satisfying $A^2\geq 0$ and $\langle \kappa_X,
A\rangle<0$ is represented by a simple $J$-holomorphic sphere (or
indeed any $J$-holomorphic curve) for any almost complex structure
$J$. To see this,  let $X_0$ be a minimal model for $X$, and let
$\pi\co X\to X_0$ be the blowdown map. Let $E_1,\ldots,E_k$ be
homology classes of the exceptional divisors of the blowups, so that
we can write $A=\pi^!A_0+\sum a_iE_i$ and
$PD(\kappa_X)=PD(\pi^*\kappa_{X_0})+\sum E_i$.  Now if our claim
were false positivity of intersections of $A$ with the $E_i$ would
show that $a_i\leq 0$, so that $\langle \kappa_X,A\rangle=\langle
\kappa_{X_0},A_0\rangle-\sum a_i\geq \langle
\kappa_{X_0},A_0\rangle$; also $A^2=A_{0}^{2}-\sum a_{i}^{2}\leq
A_{0}^{2}$.  So the assumption that $A^2\geq 0$ implies that also
$A_{0}^{2}\geq 0$, and so $A_0$ lies in the closure of the forward
light cone determined by the symplectic form on $X_0$. But Theorems
A and B of \cite{Liu} show that, since $X_0$ is minimal and neither
rational nor ruled, $\kappa_{X_0}$ is also in the closure of the
forward light cone, and so $\langle\kappa_{X_0},A_0\rangle\geq 0$,
whence $\langle \kappa_X, A\rangle\geq 0$, contradicting our
assumption on $\langle\kappa_X,A\rangle$.

Thus in this case, just as in Case 1, the sufficient condition in
Lemma \ref{adj} is vacuously satisfied. Combining this with the
results in Cases 1,2,3, and 4 now shows:

\begin{prop}\label{proveknef} Let $(X,\omega)$ be a symplectic $4$-manifold, and
$F$ an embedded symplectic surface in $F$ with genus $g>0$.  Assume
that there are no  symplectic spheres of square $-1$ contained in
$X\setminus F$, and that $X$ is not the total space of an
$S^2$-bundle over a surface of genus $g$ having $F$ as a section.
Then $F$ is rationally $K$-nef.\end{prop}

Combining Proposition \ref{proveknef} with Proposition \ref{knefmin}
thus immediately proves Case (iii) of Theorem \ref{main}.

Since Case (i) of Theorem \ref{main} is trivial, all that remains
now is to prove the statement of Case (ii).  So assume that $X_1$ is
a ruled surface with section $F_1$.  If $X_2$ is not minimal, let
$E$ be a symplectic $(-1)$-sphere in $X_2$.  If necessary, perturb
$E$ so that its intersections with $F_2$ are transverse, say at
$p_1,\ldots,p_k$.  So $E\setminus \nu_2\subset X_2\setminus \nu_2$
is a sphere with $k$ boundary components, which lie over
$p_1,\ldots,p_k$ in the circle bundle over $F_2$ which is the
boundary of $X_2\setminus \nu_2$.  Now $F_1$ contains symplectic
spheres $f_1,\ldots,f_k$ of square zero meeting $F_1\cong F_2$ at
just the points $p_1,\ldots, p_k,$ respectively, and we may (perhaps
after a small isotopy) use the discs $f_i\setminus \nu_1$ to cap off
the boundary components of $E\setminus \nu_2$ in the symplectic sum
$Z$, thus obtaining a sphere which has square $-1$ since $[E]^2=-1$
and each $[f_i]^2=0$; one can see that the sphere can be taken
symplectic by means of the pairwise sum construction from Theorem
1.4 of \cite{G}.

Conversely, if the symplectic sum $Z$ is not minimal we can simply
note that if we take the sympelctic sum of $Z$ with the ruled
surface $X_1$ by identifying the image of $F_1$ in
$Z=X_1\#_{F_1=F_2}X_2$ with $F_1\subset X_1$ using the inverse of
the gluing map that was used to form $Z$, then the resulting
manifold is deformation equivalent to $X_2$; effectively, in forming
$Z$ we have removed a neighborhood of $F_2\subset X_2$ and then
glued it back in a possibly new way, and in taking this second
symplectic sum with $X_1$ we are just regluing this neighborhood of
$F_2$ in its original configuration. So since by the previous
paragraph non-minimality is preserved under symplectic sums with
$X_1$ along $F_1$, we deduce that if $Z$ is not minimal then neither
is $X_2$. Another way of seeing this fact is to note that in the
relative Gromov--Witten theory argument in Section \ref{relgr} one
or the other of the disconnected curves representing the classes
$A_1$ and $A_2$ in equation (\ref{can}) must necessarily include a
sphere of square $-1$ as a component (for a virtual dimension count
precludes spheres of lower square than $-1$ from contributing), and
so since $X_1$ is minimal $X_2$ would have to contain a sphere of
square $-1$.
 This concludes the
proof of Theorem \ref{main}.

\section{Proof of the corollaries}
\begin{proof}[Proof of Corollary \ref{lef}] Let $f\co X\to S^2$ be a
symplectic Lefschetz fibration of positive genus $g$ on a symplectic
$4$-manifold which is relatively minimal and is not the projection
$\Sigma_g\times S^2\to S^2$, and let $F$ be a smooth fiber of $X$.
Choose $J$ generically from the set of almost complex structures on
$X$ which make $f$ a pseudoholomorphic map.  If $E\in
H_2(X;\mathbb{Z})$ is represented by an embedded symplectic
$(-1)$-sphere, then it has a $J$-holomorphic representative $S$, and
then $f|_S$ is a holomorphic map from $S^2$ to itself of degree
equal to $E\cap [F]$.  Hence $[F]\cap E>0$ unless $S$ is contained
in a fiber of $f$, but this latter possibility is forbidden by
relative minimality.

If $X$ admits a ruling with $F$ as a section, then since $F$ has
even square $X=\Sigma\times S^2$.  Again, let $J$ be an almost
complex structure making $f$ pseudoholomorphic, and as in the proof
of Proposition \ref{ruled} let $\Sigma$ denote the space of
unparametrized pseudoholomorphic curves in the class $[\{pt\}\times
S^2]$.  We then obtain an $S^2$-fibration $g\co X\to \Sigma$ having
$F$ as a section by sending $x\in X$ to the point in $\Sigma$
representing the unique $J$-holomorphic sphere on which $x$ lies.
Then $x\mapsto (g(x),f(x))$ gives an identification $X\cong
\Sigma\times S^2$ with respect to which $f$ appears as the
projection $\Sigma\times S^2\to S^2$, a contradiction.

So if $f_1,f_2$ are as in the statement of the Corollary, then their
smooth fibers $F_i$ ($i=1,2$) meet each embedded symplectic
$(-1)$-sphere in $X_i$ and are not sections of rulings on $X_i$.
Hence by Theorem \ref{main} the fiber sum of $f_1$ and $f_2$ is
minimal.
\end{proof}
\begin{proof}[Proof of Corollary \ref{Gompfex}] Gompf constructs $4$-manifolds $M_{a,b,G}$ with $\pi_1=G$, $c_{1}^{2}=a$, $c_2=b$
in Section 6 of \cite{G}; we just need to see that these are
minimal. The $M_{a,b,G}$ are obtained by fiber sums involving the
following building blocks:
\begin{itemize}
\item[(i)] A manifold $M_G$ which is spin (and hence minimal) and
satisfies $\pi_1(M_G)=G$, $c_{1}^{2}(M_G)=0$, $c_2(M)>0$.  The
summation occurs along a certain torus $T$ of square zero in $M_G$.
\item[(ii)] $P_1=\mathbb{C}P^2\#13\overline{\mathbb{C}P^2}$, which
contains a square zero, genus 2 curve $F$ which is obtained from an
irreducible quartic in $\mathbb{C}P^2$ having a single ordinary
double point by blowing up at the double point and at 12 other
points on the curve.  Note that where $S$ is the strict transform of
a line in $\mathbb{C}P^2$ which passes through the first blown-up
point and none of the others, we have $PD[F]=c_1(P_1)+PD[S]$.  If
$E$ is an embedded symplectic sphere of square $-1$ in $P_1$, we
have $[S]\cap[E]\geq 0$ (using positivity of intersections) and
$\langle c_1(P_1),[E]\rangle=1$, so $F$ meets $E$.
\item[(iii)] $P_2=\mathbb{C}P^2\#12\overline{\mathbb{C}P^2}$, which
contains a square zero, genus 2 curve $F$ which is obtained from a
sextic in $\mathbb{C}P^2$ having eight ordinary double points by
blowing up at the double points and at four other points on the
curve.  Again, since $PD[F]-c_1(P_2)$ is represented by an effective
divisor (namely the strict transform of a cubic passing through the
first eight blown-up points), we deduce as in the preceding case
that $F$ meets every embedded symplectic sphere of square $-1$ in
$P_2$.
\item[(iv)] $Q_1=(T^2\times T^2)\#2\overline{\mathbb{C}P^2}$, which
contains a square zero, genus 2 curve $F$ obtained by symplectically
resolving the double point in $\mathbb{T}^2\times\{p\}\cup
\{p\}\times T^2\subset T^2\times T^2$ and then blowing up two points
on the resulting surface.  Since $T^2\times T^2$ is aspherical the
only embedded symplectic $(-1)$-spheres in $Q_1$ are the two
exceptional divisors of the blowup, each of which meets $F$.
\item[(v)] $Q_2$, which is a torus bundle over a suface of genus 2
(hence is aspherical) and contains a square zero surface $F$ of
genus 2.
\item[(vi)] $S_{1,1}$, which is constructed as follows: let
$T_1=T^2\times\{p\},T_2=\{p\}\times T^2\subset T^2\times T^2$, and
blow up at $(p,p)$ and at 8 points on $T_1\setminus T_2$, and take
the symplectic sum with $\mathbb{C}P^2$ by identifying (the strict
transform of) $T_1$ with a cubic; the result is minimal by Theorem
\ref{main}.  Now blow up at 8 more points on the strict transform of
$T_2$ and form $S_{1,1}$ as the symplectic sum with $\mathbb{C}P^2$
by identifying the new strict transform of $T_2$ with a cubic; again
this is minimal by Theorem \ref{main}.
\item[(vii)] $\mathbb{C}P^2\#8\overline{\mathbb{C}P^2}$, which contains
a square $1$, genus 1 curve which is Poincar\'e dual to $c_1$ and
hence meets each sphere of square $-1$.
\end{itemize}
None of these building blocks are ruled, and each is either minimal
or has the property that each of its embedded symplectic
$(-1)$-spheres pass through the genus-1 or 2 surface along which the
fiber sum is taken.  Further any symplectic sum of two or more of
them along the relevant surfaces has positive Euler characteristic
and so is not a ruled surface over a curve of positive genus.  So if
$M_{a,b,G}=X_1\#_{F_1=F_2}X_2\#\cdots\# X_{k-1}\#_{F_{k-1}=F_k}X_k$
where the $X_i$ are as above, applying Theorem \ref{main}
inductively shows that for each $l\geq 2$
$X_1\#_{F_1=F_2}X_2\#\cdots\#X_{l-1}\#_{F_{l-1}=F_l}X_l$ is minimal,
and in particular $M_{a,b,G}$ is minimal.
\end{proof}

\end{document}